\documentclass{article}      
\usepackage{amssymb}
\usepackage{amsmath}
\usepackage{amscd}
\usepackage{amsthm}  
\usepackage{epsfig}  
\usepackage{subfigure}
\usepackage{floatfig}
\usepackage{wrapfig}
\newtheorem{theorem}{Theorem}
\newtheorem{corollary}[theorem]{Corollary}
\newtheorem{lemma}[theorem]{Lemma}

\newtheorem{definition}{Definition}
 
\newcommand{\bc}{\mathbb{C}}

\newcommand{\bz}{\mathbb{Z}}
\newcommand{\br}{\mathbb{R}}
\newcommand{\bq}{\mathbb{Q}}

\newcommand{\cc}{\mathcal{C}}

\newcommand{\cm}{\mathcal{M}}
\newcommand{\cf}{\mathcal{F}}

\newcommand{\hk}{\hookrightarrow}
\newcommand{\bg}{\bigskip}
\newcommand{\med}{\medskip}
\newcommand{\la}{\longrightarrow}
\newcommand{\bfl}{\begin{flushleft}}
\newcommand{\efl}{\end{flushleft}}

\newcommand{\tlm}{LM_\pm^{-TLM}}
\newcommand{\sig}{\Sigma}
\newcommand{\G}{\Gamma}
 \newcommand{\rn}{(\rho_{in})_!}
 \newcommand{\hm}{h^{mid}_*}  
\newcommand{\bph}{\bar \phi}
\newcommand{\mtm}{M^{-TM}}
\newcommand{\lpm}{LM_\pm}

\newcommand{\gro}{Gr^0_{res}}
\newcommand{\spq}{\cc \cf^\mu_{p,q}(g)}

\newcommand{\om}{\Omega}
\begin{document}  
\initfloatingfigs

\title{A polarized view of string topology}
\author{Ralph L. Cohen   \thanks{The  first author was partially supported by a grant from the NSF
}  \quad and  \quad V\'eronique Godin\\
 Stanford University 
 }
   \date{\today}
 \maketitle  
 \begin{centerline} {\sl Dedicated to Graeme Segal on the occasion of his 60th birthday.}\end{centerline}


 \begin{abstract}
Let $M$ be a closed, connected manifold, and $LM$ its loop space. 
In this paper we describe closed string topology operations in $h_*(LM)$, where $h_*$ is a generalized homology theory
that supports an orientation of  $M$.  We will show that these operations give $h_*(LM)$ the structure
of a unital, commutative Frobenius algebra without a counit. Equivalently they describe a positive boundary, two dimensional topological 
quantum field theory associated to $h_*(LM)$. This implies that there are operations corresponding to any surface with $p$ incoming and $q$ outgoing
boundary components, so long as $q \geq 1$. The absence of a counit follows from the nonexistence of an operation associated to the disk,
$D^2$, viewed as a cobordism from the circle to the empty set.  We will study   homological obstructions to constructing such an operation, and show
that in  order for such an operation to exist,  one must take
$h_*(LM)$ to  be an appropriate  homological pro-object associated to the loop space.  Motivated by this, we introduce   a prospectrum associated to
$LM$    when $M$ has an almost complex
structure.  Given such   a manifold  its loop space  has      a canonical  polarization of its tangent bundle,  which is the fundamental
feature needed to define this prospectrum. We refer to this as the ``polarized Atiyah - dual" of
$LM$.  An appropriate homology theory applied to this prospectrum would be a candidate for a theory that supports string topology operations
associated to any surface, including closed surfaces. 
\end{abstract}


\section*{Introduction}  
Let $M^n$ be a closed, oriented manifold of dimension $n$, and let $LM$ be its free loop space.  The ``string
topology" theory of Chas-Sullivan \cite{chassull}  describes  a rich structure in the homology and equivariant homology of
$LM$.   The most basic operation is an intersection - type product,
$$
\circ : H_q(LM) \times H_r(LM) \la H_{q+r-n}(LM)
$$
that is compatible with both the intersection product in the homology of the manifold, and the Pontrjagin
product in the homology of the based loop space, $H_*(\om M)$.  Moreover this product structure extends
to a  Batalin-Vilkovisky algebra structure on $H_*(LM)$, and an induced Lie algebra structure on the
equivariant homology, $H^{S^1}_*(LM)$. 
More recently Chas and Sullivan \cite{CS2} described a Lie bialgebra structure on the rational reduced equivariant homology,
$ H^{S^1}_*(LM, M;  \bq)$, where $M$ is  embedded as the constant loops in $  LM $.  

These string topology operations and their generalizations are parameterized by combinatorial data related to fat graphs used in
studying Riemann surfaces \cite{chassull},  \cite{CS2}, \cite{ sull},  \cite{vor}, \cite{cj}.  The associated field theory aspects
of string topology is a subject that is still very much under investigation.
In this paper we contribute to this investigation in the following two ways.

Recall that a two dimensional topological quantum field theory associates to an oriented compact one manifold $S$ a
vector space
$A_S$,   and to any oriented cobordism $\sig$ between  $S_1$ and
$S_2$ a linear map $\mu_{\sig} : A_{S_1} \to A_{S_2}$.   Such an assignment is required to satisfy various well known
axioms, including  a gluing axiom.  
Recall also that if
$A = A_{S^1}$,   such a TQFT structure is equivalent to a Frobenius algebra structure on $A$ 
  \cite{seg},
\cite{dijkgraaf}, \cite{abrams}.  This is a unital, commutative algebra structure, $\mu : A \otimes A \to A$, together
with   a counit (or trace map ) 
$\theta : A \to k$ so that the composition $\theta \circ \mu : A\otimes A \to A \to k$ is a nondegenerate form.  From
the  TQFT point of view, the unit in the algebra, $u : k \to A$ is the operation corresponding  to a disk $D^2$ viewed
as a cobordism from the empty set to the circle, and the counit $\theta : A \to k$ is the operation corresponding to the
disk viewed as a cobordism from the circle to the empty set.     Without the counit $\theta$ a
Frobenius algebra is equivalent to a unital,  commutative algebra $A$, together with a
cocommutative  coalgebra structure,
$\Delta : A \to A \otimes A$ without counit, where $\Delta$ is a map of $A$-modules.  From the TQFT point
of view,  a noncounital Frobenius algebra corresponds to a ``positive boundary"  TQFT,
in the sense that operations $\mu_{\sig}$ are defined only when  each component of the surface $\sig$ has a \sl positive \rm number
of outgoing boundary components.   
 
   Let $h^*$ be a multiplicative generalized cohomology theory whose coefficient ring, $h^*(point)$ is a  graded field (that is, every nonzero homogeneous
element is invertible).  Besides usual cohomology with field coefficients, natural examples of such theories are periodic $K$-theory with
field coefficients, or any of the Morava
$K$-theories.  Any  multiplicative generalized cohomology theory gives rise to such a theory by appropriately localizing
the coefficient ring.    

 Let   
$h_*$ be the associated generalized homology theory. Let $M^n$ be a  closed, $n$ - dimensional manifold  which is oriented with
respect to this theory.  Our first result, which builds on work of Sullivan \cite{sull},  is that  string topology operations can be defined  to give a 
two dimensional positive boundary TQFT, with $A_{S^1}  = h_*(LM)$.     

\begin{theorem}  \label{one}
The homology of the free loop space
$h_*(LM)$ has the structure  of a Frobenius algebra  without counit.   The ground field of this algebra structure
is the coefficient field, $h_* = h_*(point)$.
 \end{theorem}

The construction of the TQFT operations corresponding to a surface $\sig$   will involve
studying spaces of  maps from a fat graph $\G_\sig$ associated to the surface to $M$,  and viewing that space as a finite codimension
submanifold of a
$(LM)^p$, where $p$ is the number of incoming boundary components of $\sig$.   We will show that this allows the construction of a 
Thom collapse map for this  embedding, which will in turn define a push-forward map $\iota_! : h_*(LM)^{\otimes p} \to h_*(Map(\G_\sig; M))$.
 The operation $\mu_\sig$  will then be defined as the composition $\rho_{out} \circ \iota_! : h_*(LM)^{\otimes p} \to h_*(Map(\G_\sig; M) \to
h_*(LM)^{\otimes q}$ where $\rho_{out}$ is induced by restricting a map from $\G_\sig$   to its outgoing boundary components.

The second goal  of this paper is to investigate the obstructions to constructing  a
homological theory applied to the loop space  which supports the string topology operations,  \emph{and} permits the  definition of a
counit in the Frobenius algebra structure, or equivalently, would   eliminate the ``positive boundary" requirement
in the  TQFT structure.   Let $\hm (LM)$ be such a conjectural theory.  In some sense this would represent a ``middle dimensional", or
``semi-infinite" homology theory associated to the loop space,  because of the existence of a nonsingular form
$\hm (LM) \otimes \hm (LM) \to k$ analogous to the interesection form on the middle dimensional homology of an even dimensional oriented manifold.

We will
see that  defining a counit would involve the construction of a push-forward map for the embedding of constant loops $M \hk LM$.  Unlike the embeddings
described above, this  has
\emph{infinite} \rm codimension.   We will argue that this  infinite dimensionality will force the use of an inverse limit of homology theories, or a
``pro-homology theory" associated to the loop space.   Using previous work of the first author and Stacey \cite{cs}, we will show that there are obstructions
to the construction of  such a pro-object unless $M$ has an almost complex structure.  In this  case the tangent bundle of the loop space has a canonical
complex polarization, and we will use it to define       the ``polarized loop space",  $LM_\pm$.  This space
  fibers over $LM$,  where an element in the fiber over
$\gamma \in LM$ is a representative of the polarization of the tangent space $ T_\gamma LM$.   We will examine various properties of
$LM_\pm$, including its equivariant properties.  We will show that the pullback of the tangent bundle  $TLM$ over $LM_\pm$ has a filtration that will allow us to
define  
  a prospectrum
$\tlm$, which  we call   the  (polarized) ``Atiyah dual" of $LM$.  We will end by describing  how the application of an appropriate    
equivariant homology theory to this prospectrum should be a good candidate for studying further field theory properties of string
topology.  This will be the topic of future work, which will be joint with J. Morava and G. Segal.

The paper is organized as follows.  In section 1 we will describe the type of fat graphs needed to define the string topology operations.  These are chord
diagrams of the sort introduced by Sullivan \cite{sull}.  We will define the topology of these chord diagrams using categorical ideas of Igusa
\cite{igusa1}\cite{igusa2}.  Our main technical result, which we will need to prove the invariance of the operations, is that the space of chord diagrams
representing surfaces of a particular diffeomorphism type is connected.   In section 2 we   define the string topology operations
and prove theorem 1.  The operations will be defined using a homotopy theoretic construction (the ``Thom collapse map")  generalizing what was done by the first
author and Jones in
\cite{cj}.   In section 3 we  describe the obstructions to the existence of a counit or trace in the
Frobenius algebra structure.    Motivated by these observations, in section 4 we describe  the ``polarized Atiyah dual" of  the loop space of an almost complex 
manifold, and study its properties. 

The authors are grateful to J. Morava, G. Segal, and D. Sullivan for many inspiring conversations about the topics of this paper.  They are also grateful to
the referee for many useful suggestions.

 \section{Fat graphs and Sullivan chord diagrams}
 Recall from  \cite{penner}
\cite{strebel} that a   fat graph    is a graph whose vertices are at least trivalent, and where the edges coming into each vertex come
equipped with a cyclic ordering.    Spaces of  fat graphs have been used   
 by many different authors  as an extremely effective tool in studying the topology and geometry of moduli spaces of Riemann surfaces.  
The essential feature of a fat graph is that when it is  thickened, it produces a   surface with boundary, which is well defined up to homeomorphism.

For our purposes, the most convenient approach to the space of fat graphs is the categorical one described by Igusa   
\cite{igusa1},\cite{igusa2}.  In  \cite{igusa1}  (chapter 8) he defined a   category  $\cf at_n(g)$    as follows.  
The objects  of $\cf at_n(g)$ are fat graphs (with no lengths assigned to the edges) , and  the  morphisms are  maps of fat graphs $f : \G_1 \to
\G_2$ (i.e maps of the underlying simplicial complexes that preserve the cyclic orderings) satisfying the following properties:
\newcounter{point}
\setcounter{point}{1}
\begin{list}{(\alph{point})}{\setlength{\itemsep}{0pt}}
\item The inverse image of any vertex is a tree.
\refstepcounter{point}
\item The inverse image of an open edge is an open edge.
\end{list}

Igusa proved that the geometric realization   $|\cf at_n (g)|$ is homotopy equivalent  to the classifying
space, 
$B\cm_{g,n}$, where
$\cm_{g,n}$ is the mapping class group of genus $g$ surfaces with $n$ marked, ordered points (theorem  8.6.3 of \cite{igusa1}).  In this
theorem $n \geq 1$ for $g \geq 1$ and $n \geq 3$  for $g=0$.  He also proved (theorem 8.1.17) that $|\cf
at_n (g)|$  is homotopy equivalent to the space of
\sl metric fat graphs,
\rm  which we denote
$\cf _n(g)$,  which is a simplicial space made up of fat graphs with appropriate metrics.   These  spaces are closely related the simplicial
sets  studied by Culler and Vogtmann
\cite{culvogt} and Kontsevich
\cite{kont}.  See
\cite{igusa1} chapter 8 for details.

 Following \cite{culvogt}  there are ``boundary cycles" associated to a fat graph $\G$ defined as follows.    Pick an edge and choose an orientation on it. 
Traversing that edge in the direction of its orientation leads   to a vertex.   Proceed with the next edge    emanating from
that vertex in the cyclic ordering, and give it the orientation pointing away from that vertex.  Continuing in this way, one
traverses several edges, eventually returning  to the original edge, with the original orientation.  This yields a ``cycle" in the set of
oriented edges and represents a boundary component.  One
partitions the set of all oriented edges into cycles, which enumerate the boundary components of the surface represented by 
$\G$.  The cycle structure of the oriented edges completely determines the combinatorial data of the fat graph. 

In a metric fat graph each boundary cycle has an orientation and a metric on it.    Hence introducing
a marked point for each boundary cycle would yield a parameterization of the boundary components. Notice that 
it is possible that two marked points    lie on the same edge, and indeed a single point on an edge might have a
``double marking" since a single edge with its two orientations might lie in two different boundary cycles.  

We call the space of metric fat graphs representing surfaces of genus $g$ with $n$ boundary components, that come equipped with
marked points on the the boundary cycles,
$\cf^\mu_n ( g)$.   This is the space of \emph{marked} metric fat graphs.  Using Igusa's simplicial set construction, one
sees that $\cf^\mu_n ( g)$ can be topologized so that the projection map that forgets the markings, 
\begin{equation}\label{forget}
p :  \cf^\mu_n ( g) \to  \cf_n (g)
\end{equation}
is a quasifibration whose fiber is the space of markings on a fixed fat graph, which is homeomorphic to  the torus $(S^1)^n$. 
The topology of the space of marked metric fat graphs is studied in detail in \cite{vero} with applications to specific combinatorial
calculations.

\med
For the purposes of constructing the string topology operations, we will use a particular type of fat graph due to Sullivan.

\med
\begin{definition}   A  ``Sullivan chord diagram"  of type $(g; p,q)$ is a  fat graph  representing a surface of genus $g$ with
$p+q$ boundary components, that consists of  a disjoint union of 
$p$  disjoint closed circles   together with the disjoint union of connected trees
whose endpoints lie on the circles.    The cyclic orderings  of the edges at the vertices must be such that each of the $p$ disjoint
circles is a boundary cycle.   These $p$ circles are referred to as the incoming boundary cycles, and the other $q$ boundary cycles
are referred to as the outgoing boundary cycles. 
\end{definition}

The ordering at the vertices in the diagrams that follow
are indicated by the  clockwise cyclic ordering of the plane.  
Also in a Sullivan chord diagram, the vertices and edges that lie on one of the $p$ disjoint circles will
be referred to as circular vertices and circular edges respectively.  The others will be referred
to as ghost vertices and edges.   
 
\begin{figure}
\begin{center}
\epsfig{file=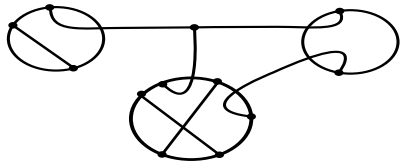, height=60pt, width=150pt}\caption{Sullivan chord diagram of type (1;3,3)}
\end{center}
\end{figure}

To define the topology on the space of metric chord diagrams, we first need to define the space of metric fat graphs, $\cf_{p,q}(g)$, of genus $g$, with
$p+q$ ordered boundary cycles,   
   with the first $p$  distinguished as incoming, and the remaining $q$ distinguished as
  outgoing.    
Igusa's simplicial construction of $\cf_n(g)$ defines a  model of $\cf_{p,q}(g)$ as the geometric realization of a simplicial set.   Moreover this space is
homotopy equivalent to the realization of the nerve of the category $\cf at_{p,q}(g)$ defined as above, with the additional feature that the objects come
equipped with an ordering of the boundary cycles, with the first
$p$ distinguished as incoming cycles.    The morphisms must preserve this structure.   

Consider the space of ``metric chord diagrams", $\cc \cf_{p,q}(g)$ defined to be the subspace of the metric fat graphs $\cf_{p,q}(g)$ whose underlying
graphs are Sullivan   chord diagrams  of type $(g; p,q)$.  
So if  $\cc \cf at_{p,q}(g)$ is the full subcategory of $\cf at_{p,q}(g)$ whose objects are chord diagrams of type $(g; p,q)$, then Igusa's argument shows
that the space of metric chord diagrams $\cc \cf_{p,q}(g)$ is homotopy equivalent to the realization of the nerve of the category
$\cc\cf at_{p,q}(g)$.

\begin{figure}  
\begin{center}
\mbox{
\subfigure[c]{\epsfig{file=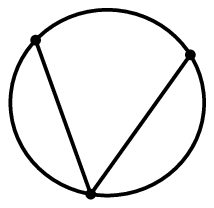, height=60pt, width=60pt}}\qquad
\subfigure[S(c)]{\epsfig{file=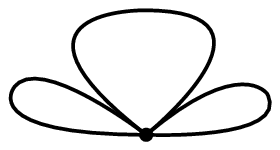, height=60pt, width=70pt}}}
\caption{Collapsing of ghost edges}\label{fig:S(c)}
\end{center}
\end{figure}

Given a metric chord diagram $c \in  \cc \cf_{p,q}(g)$, there is an associated metric fat graph, $S(c)$, obtained from $c$ by collapsing all ghost edges.  
There is an  induced cyclic ordering on the vertices of $S(c)$ so that the collapse map $\pi : c \to S(c) $ is a map of fat graphs in $\cf
at_{p,q}(g)$. 
 Figure
\ref{fig:S(c)} describes  this collapse map.
 
We will define a  marking of a Sullivan chord diagram $c$ to be a marking (i.e a choice of point) on each  of the boundary cycles of  
the associated fat graph $S(c)$. 
We let $\cc \cf^\mu_{p,q}(g)$ denote the space of all marked metric chord diagrams.  Like with the full space of marked fat graphs, this space can be
topologized in a natural way so that the projection map that forgets the markings,
\begin{equation}\label{unmark}
p : \cc \cf^\mu_{p,q}(g) \la \cc \cf_{p,q}(g)
\end{equation}
is a quasifibration, with fiber over a metric chord diagram $c$ equivalent to a torus $(S^1)^{p+q}$.  
Again, the topology of these spaces of marked chord diagrams will be studied in detail in \cite{vero}.

The space of marked, metric chord diagrams $\spq$ will be used in the next section to parameterize the string topology operations.   Its topology,  
  however,  is far from understood.  It is a proper subspace of a space homotopy equivalent to the classifying space of the mapping class group, and
thus moduli space.  However very little is known about the topology of this subspace.  We make the following conjecture, which would say that the
parameterizing spaces of string topology operations are homotopy equivalent to moduli spaces of curves, thus potentially leading to a  
conformal field theory type structure.

\med
 {\bf Conjecture.} The inclusion   $\spq \hk  \cc \cf_{p,q}(g)$  is a homotopy equivalence, and in particular 
$\spq$ is homotopy equivalent to the classifying space
$B\cm_{g}^{p+q}$, where $\cm_g^{p+q}$ is the mapping class group of a surface of genus $g$ with $p+q$
ordered boundary components, and the diffeomorphisms preserve the boundary components pointwise. 
  
\med
For the purposes of this paper we will need the following property of these spaces.

\begin{theorem}\label{connect}
 The space $\spq$ is path connected.
\end{theorem}
 
\emph{Proof.}    Because of   quasifibration (\ref{unmark}), it suffices to prove that the  space of \emph{unmarked} metric chord
diagrams
$\cc \cf_{p,q}(g)$ is connected.   However as remarked earlier,  this space is homotopy equivalent 
to the nerve of the category $\cc \cf at_{p,q}(g)$.

\begin{wrapfigure}{r}{210pt}
\epsfig{file=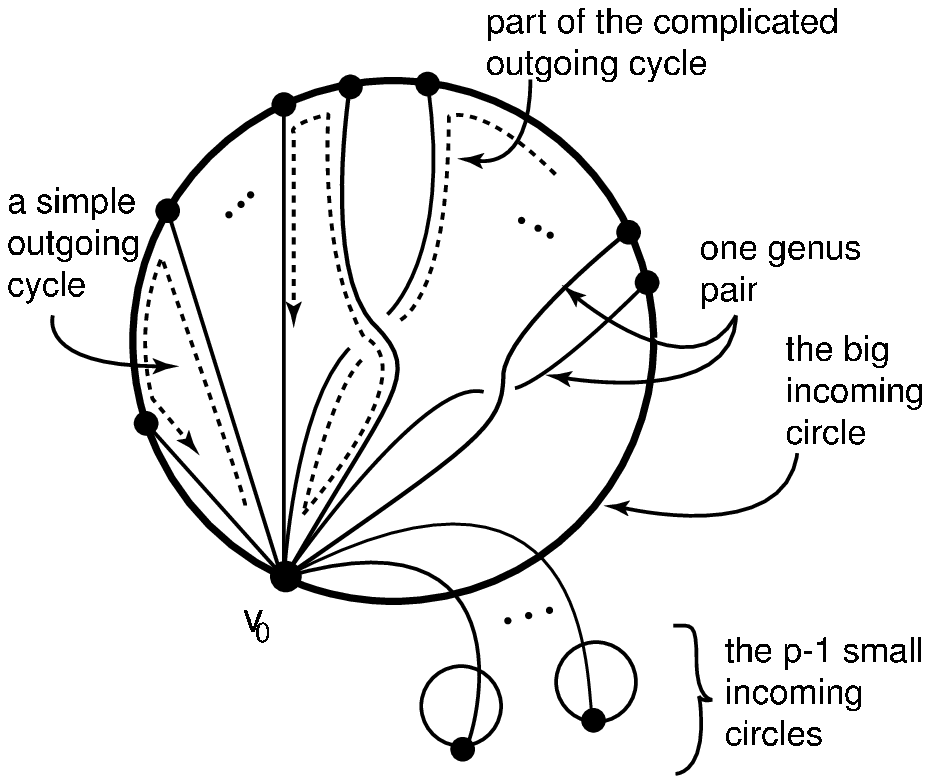, height=150pt, width=200pt}\caption{$\Gamma_0$}\label{fig:basepoint}
\end{wrapfigure}
Now a morphism  between objects $a$ and $b$ in a category determines a path from $a$ to $b$ in the   
geometric realization of its nerve. Since a morphism in  $\cc \cf at_{p,q}(g)$ collapses trees to vertices, we refer to such morphisms as ``collapses".  Now
by   reversing orientation, a morphism from $b$ to $a$ also determines a path from $a$ to $b$ in the geometric realization.  We refer to such a   
    morphism as an ``expansion"  from $a$ to $b$.    Therefore   to prove this theorem  it
suffices to build, for any chord diagram in our space, a sequence of collapses and  and expansions from it to a fixed chord diagram 
$\Gamma_0$.    In the following diagrams, dashed lines will represent boundary cycles.

We will choose our basepoint $\Gamma_0$ as in figure \ref{fig:basepoint}.  In $\Gamma_0$, $p-1$ of the
incoming circles contain  only one vertex.  There is also a distinguished vertex $v_o$ in the $p^{th}$ incoming circle (which we refer to as the
``big circle"). Moreover 
$q-1$ outgoing components share the same structure :  they can be traced by going from $v_o$ along a chord edge whose other vertex also lies
on the big circle, going along the next circle edge in the cyclic ordering, and then going along the next chord edge back to $v_o$. 
In the last outgoing boundary component positive genus is produced by  pairs of chord edges twisted, as shown. (These pairs add two
generators in the fundamental group of the surface but do not affect the number of boundary components, therefore they ``create
genus".) Notice also that except for $v_o$ all of the vertices of $\Gamma_0$  are trivalent.
So, in $\Gamma_0$ the complexity is concentrated in the big incoming circle, the last outgoing boundary circle,  and one vertex $v_o$.

The ordering of the boundary cycles in $\Gamma_0$  is given by making the first incoming cycle the one containing $v_o$.   The 
ordering of the other incoming boundary components follows the cyclic ordering at $v_o$
(so that in figure \ref{fig:basepoint} the circle on the right will come second and the one on the left last.) Similarly the cyclic ordering
at $v_o$ will give us an ordering of the outgoing boundary components (in which the complicated boundary cycle is last). 

To prove  the theorem we  start from any chord diagram (object) in  $\cc \cf at_{p,q}(g) $ and get to   $\Gamma_0$ by a sequence of collapses and
expansions.  In our figures, the arrows follow the direction of the corresponding morphism in our category $\cc \cf at_{p,q}(g)$.
Note that since, in a Sullivan chord diagram, the incoming boundaries are represented by disjoint circles, a chord edge between two
circular vertices cannot be collapsed.  Remember also that the ghost edges need to form a disjoint union of trees. Hence if both vertices 
of a circular edge are part of the same tree of ghost edges (same connected component of the ghost structure), this circular edge cannot be collapsed. We
will call an edge ``essential'' if it cannot be collapsed. That is, it  is either a  circular edge and its collapse would create a non-trivial
cycle among    the ghost edges, or it is a chord edge between two circular vertices. 

Throughout this proof, letters from the beginning of the alphabet will be used to  label edges that are on the verge of 
being collapsed,  and  letters from the end of the alphabet will be used to label edges that have just been created, via an expansion.
We will start by assuming that all nonessential edges have been collapsed.  


\bg
\begin{wrapfigure}{l}{190pt}
\epsfig{file=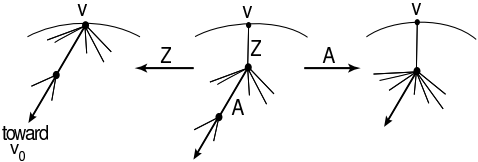, height=60pt, width=180pt}\caption{Pushing edges}\label{fig:gamma}
\end{wrapfigure}
The first step will be to find a path to a chord diagram with a  distinguished vertex $v_o$, 
the only one with more than three  edges emanating from it.  Choose $v_o$ to be any 
vertex on the first incoming boundary cycle. 
For any vertex $v$, other than $v_o$, having more than three edges,    we will  ``push"  the edges of  $v$  toward  $v_o$ by
a sequence of expansions and collapses.  This is done as follows. 

Since all edges are essential, the vertices of any circular edge are part of the same connected
component of the ghost structure. 
 We can therefore choose a path $\gamma$ from $v$ to $v_o$ contained completely in the the ghost structure.  Following figure
\ref{fig:gamma}, we can push the edges of $v$ a step closer to $v_o$.  
 Repeating this process completes this step.

We now have a distinguished vertex $v_o$, which is  the only vertex with more than one ghost edge. Note that there is a unique ordering of the edges
emanating  from $v_o$ that is compatible with the cyclic order such that one of the circular edges of $v_o$ 
is first in the order, and the other circular edge is last in the order. We will think of this ordering of edges as passing from left to right.


\bg
\begin{wrapfigure}{r}{190pt}
\epsfig{file=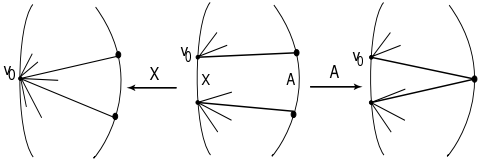, height=60pt, width=180pt}\caption{Getting rid of one edge}\label{fig:1incoming}
\end{wrapfigure}
We will next simplify the incoming circles.  In this step all of the incoming cycles but the first will be brought down to one edge. 
Notice first that all ghost edges have $v_o$ as one of their vertices.
A ghost edge between two other vertices would be in itself a connected component of the ghost structure.
But a circular edge joining two different components can be collapsed without creating a ``ghost cycle''.
Since all the edges of our graphs are essential, we know that all ghost edges have $v_o$ as a vertex.

Now take any ``small'' incoming boundary circle containing more than one edge . 
As seen in figure \ref{fig:1incoming}, the addition of an edge $X$ close to
$v_o$, renders $A$ nonessential and it can be collapsed.

\begin{figure}
\begin{center}
\epsfig{file=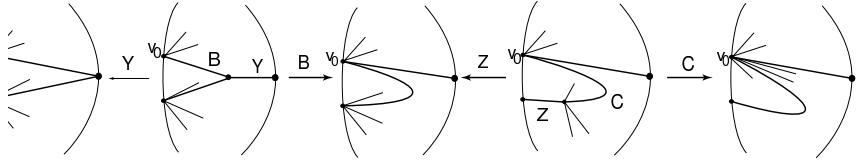, height=60pt, width=300pt}\caption{Pushing edges towards $v_o$}\label{fig:2incoming}
\end{center}
\end{figure}

After this has been done, there will be three non-trivalent vertices.
The procedure shown in figure \ref{fig:2incoming} brings this number back down to one.   
Observe that in this procedure there is no risk of collapsing an essential edge. Repeating this process
reduces the number of vertices on these incoming circles down to one per circle. Notice that our chord
diagram still has a unique non-trivalent vertex $v_o$.

\begin{figure}
\begin{center}
\epsfig{file=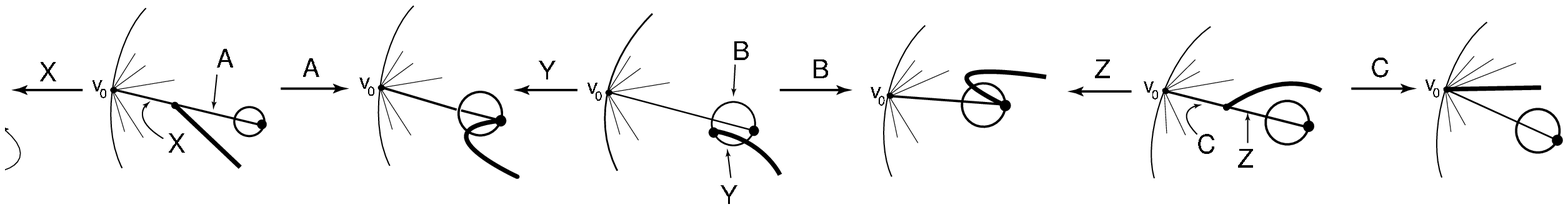, height=60pt, width=400pt}\caption{Switching an incoming circle and an edge}\label{fig:switch}
\end{center}
\end{figure}

Now the $p-1$ small incoming circles have a unique vertex with a unique ghost edge linking it with $v_o$ on the big incoming circle. 
Using their simple structure we will be able to switch the order at $v_o$ of their ghost edges and any other ghost edge
(see figure \ref{fig:switch}).  


\bg
\begin{wrapfigure}{l}{100pt}
\epsfig{file=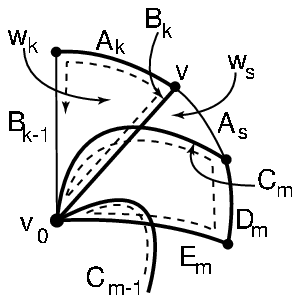, height=90pt, width=90pt}\caption{}\label{fig:kth}
\end{wrapfigure}
The next step is to simplify the structure of $q-1$ outgoing boundary components (all but the last one). This will be achieved by lowering the
number of edges involved in the tracing of these boundary cycles to a minimum (three edges for
most).  We will use the term ``clean" to refer to this simple form.   
  
Now
each of the outgoing boundaries has an edge on the big circle. (The cleaned  boundaries will only have one such edge.)
Assume by induction that the first $k-1$ outgoing boundaries have already been cleaned. 
Assume also that these cleaned boundaries have been pushed to the left of the big incoming circle, meaning that
unique incoming edge is situated to the left of all the incoming edges associated to the uncleaned boundaries (no uniqueness here) and that their
ghost edges at $v_o$ are attached left of all other ghost edges. 
Assume that the cleaned boundaries have been ordered.
Assume also that at least two outgoing boundaries still need to be cleaned. 
Firstly we will clean the next outgoing boundary $w_k$ and secondly we will push it left to its proper position.

Since all the clean boundaries are to the left, we will have two successive circular edges, $A_k$ and $A_s$, 
on our big circle such that $A_k$ is part of $w_k$ the next outgoing boundary component to be cleaned,
and $A_s$ belongs to $w_s$ with $s>k$. We will argue the case where $A_k$ is to the left of $A_s$. The other situation is argued similarly.

\begin{figure}
\begin{center}
\epsfig{file=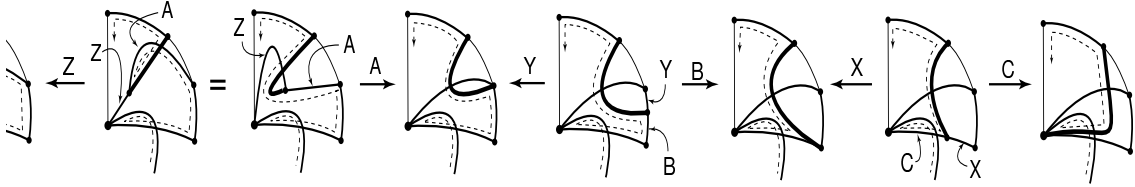, height=70pt, width=420pt}\caption{Skipping three edges}\label{fig:genus}
\end{center}
\end{figure}

Let $v$ be the common circular vertex of $A_k$ and $A_m$ and let $B_k$ be the ghost edge coming into $v$.
Notice that  the cycles representing both $w_k$ and $w_m$ include $B_k$ 
(with    different orientations).
We   now push all the extra edges involved in the tracing of $w_k$ to the right of $B_k$ and hence into $w_m$. 
Since all the ghost edges involved in the cycle representing $w_k$ start at $v_o$ and end on the big circle, 
the cycle traced by $w_k$ is ($B_{k-1}$,...,$C_{m-1}$, $E_m$, $D_m$, $C_m$, $B_k$, $A_k$) where the $E$'s 
and the $C$'s are ghost edges from $v_o$ to the big circle, and the $D$'s are circular edges on
the big circle.  See figure \ref{fig:kth}. In figure \ref{fig:genus} the start of $B_k$ is glided
along the edges $C_m$ $D_m$ and $E_m$. ($B_k$ is thickened on each of the pictures to help
visualize this process.) After these glidings $w_k$ will not include $E_m$, $D_m$ and $C_m$.
We can repeat this process until only $B_k$, $A_k$ and $B_{k-1}$ are left in the tracing of
$w_k$.

\begin{figure}
\begin{center}
\epsfig{file=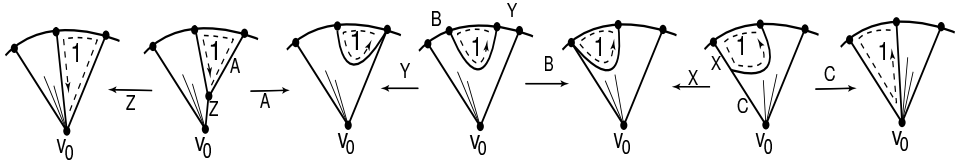, height=60pt, width=360pt}\caption{Moving boundary 1 to the left}\label{fig:moving}
\end{center}
\end{figure}
To make sure we don't interfere with this boundary cycle when making subsequent rearrangements, we will move it to the left of all the remaining 
uncleaned boundaries. 
Follow figure \ref{fig:moving} to see how to switch this boundary (labelled 1) with the one directly on the left of it. 
By induction we will can clean all of the outgoing boundary but the last one.

\begin{wrapfigure}{l}{200pt}
\epsfig{file=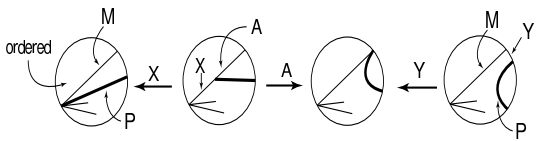, height=50pt, width=200pt}\caption{Moving $P$ away}\label{fig:untangle1}
\end{wrapfigure}

We are now very close to our goal. We have the incoming and the outgoing boundaries in the right order and in the right form. 
The only issue rests with the last outgoing boundary component. The one that includes all the ``genus creating'' edges. 
To finally reach $\Gamma_0$, we need to untangle these into twisted pairs . This is done by induction on the number of pairs of such edges 
left to untangle. 

Choose $M$ to be first genus creating edge (in the ordering at $v_o$) that has not yet been paired in our inductive process. Our first 
step is to find an edge that ``crosses it'', meaning that it starts on the right of $M$ at $v_o$ and ends up on the left of $M$ on the big circle.
Take $P$ to be the next edge at $v_o$. 
If $P$ ends up on the left of $M$, we have our edge and we are ready to apply the second step.
If this is not the case, we'll move $P$ along $M$ as shown in figure \ref{fig:untangle1} and consider the next edge.
Since both orientation of $M$ are part of the last boundary component, the tracing of this component
moves from the edges on the right of $M$ to the edges on the left of $M$. This implies that there is
at least one edge that starts on the right of $M$ at $v_o$ and ends on the left of $M$ on the big circle. 
So by moving through the edges $P$ we will find one of these crossing edges and this step will be completed.

\begin{figure}
\begin{center}
\epsfig{file=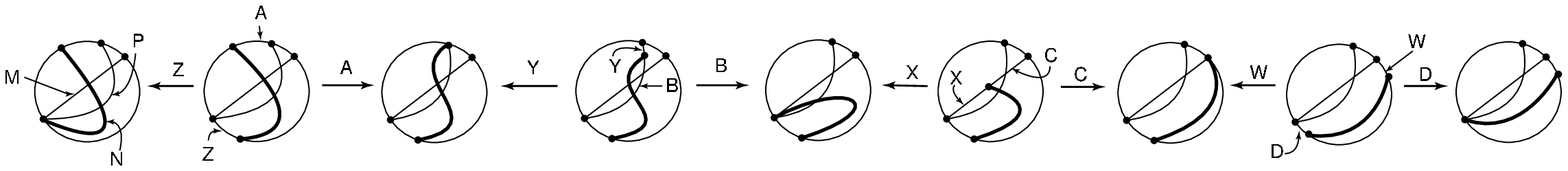, height=60pt, width=400pt}\caption{Separating the pair $M$-$P$ from $N$}\label{fig:untangle2}
\end{center}
\end{figure}

Now we have our pair of edges $M$ and $P$. But we   would like to have $M$ and $P$ completely separate from the other edges as in 
$\Gamma_0$. Any edge $N$ that is intertwined with $M$ and $P$ needs to be moved. 
Since everything on the left of $M$ and $P$ at $v_o$ is already in proper
order, this $N$ will always start on the right of $M$ and $P$ at $v_o$ (and end up on the big circle). 
But there are still two ways that $N$ might be intertwined with our pair: 
it could either end up between $M$ and $P$ or on the left of these edges.
Figure \ref{fig:untangle2} shows how to glide the edge $N$ along first $M$ and then $P$. 
A similar operation would get rid of the edges landing between $M$ and $P$. 

Before restarting these steps for isolating the next pair, we need to bring a lot of edges back to $v_o$. For example in figure \ref{fig:untangle1} $P$ 
ends up completely disconnected from $v_o$ and from the rest of the ghost structure. 
To achieve this we will first collapse all nonessential edges and then we'll reduce the number of edges on the non-$v_o$ 
vertices down to three as done in figure \ref{fig:gamma}. Note that $M$ and $P$ can be kept isolated while all of this is done.
 
\bg
After this process, all the genus edges are paired and twisted properly. We can finally put the ghost edges connecting the small circles to
$v_o$ in their correct position. 
Our chord diagram has now one big incoming circle and a special vertex $v_o$,  the only
non-trivalent vertex. It has $p-1$ small ``one-vertexed''  incoming circles linked with $v_o$ by
one ghost edge landing in the last outgoing boundary cycle and ordered properly.  It has $q-1$
simple  outgoing boundary components positioned on the left of the big incoming circle. The last
outgoing component contains all the genus edges isolated into twisted pairs and the ghost edges
linking the big incoming circles with the small ones. This means that the cycles associated to the
different boundary components are exactly the same in this chord diagram and in $\Gamma_0$.
But we know that these cycles determine completely a fat graph. We  have shown how to connect
a random chord diagram to $\Gamma_0$ by a sequence of collapses and expansions.  This proves
the theorem.    

\vspace{-15pt}
\begin{flushright} 
{$\square$}
\end{flushright}

\section{ The Thom collapse map  and string topology  operations}   

In this section we use fat graphs to define the string operations, and will    prove theorem \ref{one} stated in the introduction.

Let $LM$ denote the space of piecewise smooth maps, $\gamma : S^1 \to M$.    
Let $\sig$ be an oriented  surface of genus $g$ with $p+q$ boundary components.

For $c
\in\spq $, let $Map_*(c, M)$ be the space of continuous maps $f : c \to M$, smooth on each edge, which is constant on each ghost edge.
Equivalently, this is the full space of maps   $S(c) \to M$ where $S(c)$ is the marked  metric fat graph described in the previous section (obtained  from
$c$ by collapsing  each  ghost  component to a point).  
Since each ghost component is a tree  and therefore contractible, this mapping space is homotopy equivalent to the space of all continuous maps,
$c \to M$, which in turn is homotopy equivalent to the smooth mapping space, $Map (\sig, M)$.  Furthermore,  the markings on $S(c)$ induce   
parameterizations of  the incoming and outgoing boundary  cycles   of $c$,   so restriction to these boundary cycles induces a diagram,  

\begin{equation}\label{cbasic}
\begin{CD}
(LM)^q @<\rho_{out}<<  Map_*(c , M) @>\rho_{in} >> (LM)^p.
\end{CD}
\end{equation}

Since $Map_*(c, M)$ is the same as the space of all continuous maps $Map (S(c), M)$,   it is
clear that 
the restriction to the incoming boundary components, 
$$
\rho_{in} : Map_*(c, M) \la (LM)^p
$$
is an \sl embedding \rm of infinite dimensional manifolds, but it has finite codimension.  We now consider its normal bundle.

Let $v(c)$ be the collection of circular vertices of a chord diagram $c$.  Let $\sigma (c)$ be the collection of vertices of the associated graph
$S(c)$.  The projection map $\pi : c \to S(c)$ determines a surjective set map, $\pi_*: v(c) \to \sigma (c)$.  For a vertex $v \in \sigma (c)$, we define the
\sl multiplicity, \rm $ \mu (v)$, to be the cardinality of the preimage, $\#\pi^{-1}(v)$.  Let
$M^{\sigma (c)}$ and
$M^{v(c)}$ be the induced mapping spaces from these vertex collections.  Then  $\pi$ induces  a diagonal map
$$
\Delta_c : M^{\sigma (c)} \la M^{v(c)}.
$$
The normal bundle of this diagonal embedding is the product bundle,
$$
\nu (\Delta_c) = \prod_{v \in \sigma (c)} (\mu (v)-1) TM \to \prod_{v \in \sigma (c)} M = M^{\sigma (c)}.
$$ 
Here $k\cdot TM$ is the $k$-fold Whitney sum of the tangent bundle.  Since $\sum_{v\in \sigma (c)} \mu (v) = v(c)$,  the fiber dimension of this bundle
is $(v(c) - \sigma (c))n$.  An easy exercise verifies that $(v(c) - \sigma (c)) = -\chi (\Sigma_c)$,  minus the Euler characteristic of a surface
represented by $c$.

Now remember that the markings of the    incoming boundary components of $S(c)$  define parameterizations of the incoming boundary components
of $c$, since these cycles only consist of circular edges. 
Now using these parameterizations    we can identify $(LM)^p$  with  $Map (c_1 \sqcup \cdots \sqcup c_p ,
M)$, where $c_1, \cdots , c_p$ are the $p$ incoming boundary cycles of $c$.  Consider the evaluation map  $ e_c  : (LM)^p \to M^{v(c)}$  defined on an
element of  
$\gamma
\in (LM)^p$ by evaluating
$\gamma$ on the circular vertices.   Similarly, define
$$\quad e_c : Map_*(c, M) \to M^{\sigma (c)}$$  by evaluating a map $f : S(c) \to M$ on the vertices.  These evaluation maps are fibrations, and notice
that the following is a pull-back square:

\begin{equation}\label{pullback}
 \begin{CD}  
Map_*(c, M)   @>\rho_{in} >\hk > (LM)^p \\
@Ve_c VV   @VVe_c V \\
M^{\sigma (c)} @>\hk >\Delta_c >  M^{v (c)}
\end{CD}
\end{equation}

By taking the inverse image of a tubular neighborhood of the embedding $\Delta_c$, one has the following consequence.

\med
\begin{lemma}\label{nbd}
  $\rho_{in} : Map_*(c, M)  \hk (LM)^p $ is a codimension $-\chi(\sig_c)n$  embedding, and has an open neighborhood
 $\nu (c)$  diffeomorphic to the total space of
the  pullback bundle, $e_c^*(\nu (\Delta_c)) = e_c^*(\prod_{v \in \sigma (c)} (\mu (v)-1) TM )$.  The fiber of
this bundle over a map $f : c \to M$ is therefore given by
$$
\nu (c)_f = \bigoplus_{v \in \sigma (c)} \bigoplus_{(\mu (v) -1)} T_{f(v)}M,
$$ where $\bigoplus_{(\mu (v) -1)} T_{f(v)}M$ refers to taking the direct sum  of $\mu (v) -1 $ copies of  $T_{f(v)}M$.   
\end{lemma}

\med
Let $Map_*(c, M)^{\nu( c)}$ be the Thom space of this normal bundle.  
This result allows us to define a Thom collapse map $\tau : (LM)^p \to Map_*(c, M)^{\nu( c)}$ defined, as usual, to be the identity
inside the tubular neighborhood, and the basepoint outside the tubular neighborhood. 

Now let $h^*$ be a generalized cohomology theory  as before.  By the above description of the bundle $\nu( c)$ we see that since
$M$ is $h^*$-oriented, the bundle $\nu ( c)$ is $h^*$-oriented.   This defines a Thom isomorphism,
$$
  t : h_*(Map_*(c, M)^{\nu( c)}) \cong  h_{*+\chi (\Sigma_c)n}(Map_*(c, M)).
$$

Now since we are assuming that the coefficient ring $h_* = h_*(point)$ is a graded field, the Kunneth spectral sequence collapses,
and hence
$$h_*(X\times Y) \cong h_*(X) \otimes_{h_*}h_*(Y).$$  From now on we take all tensor products to be over the ground field $h_*$. We can therefore
make the following definitions.

\begin{definition} Fix $c \in \spq.$

a.  Define the push-forward map $\rn : h_*(LM)^{\otimes p}  \to h_{*+\chi (\sig_c)n}(Map_*(c, M))$ to be the composition
$$
\begin{CD}
\rn : h_*(LM)^{\otimes p} \cong h_*((LM)^p)  @>\tau_*>>  h_ * (Map_*(c, M)^{\nu(\Delta_c)})  @>t>\cong>  h_{*+\chi (\Sigma_c)n}(Map_*(c, M)).
\end{CD}
$$

b.  Define the operation $\mu_c : h_*(LM)^{\otimes p}   \to h_*(LM)^{\otimes q} $ to be the composition,
$$
\begin{CD} 
\mu_c :  h_*(LM)^{\otimes p}  @>\rn >>   h_{*+\chi (\sig_c)n}(Map_*(c, M)) @>(\rho_{out})_*>> h_{*+\chi (\sig_c)n}((LM)^q)   .
\end{CD}
$$
\end{definition}

\med
In order to use these operations to prove theorem  \ref{one}  we will first need to verify the following.

\med
\begin{theorem}\label{invariant} The operations $\mu_c : h_*(LM)^{\otimes p}   \to h_*(LM)^{\otimes q} $  do not depend on the choice of marked metric
chord diagram
 $c  \in  \spq$.  In other words, they only depend on the topological  type (g; p,q) of the chord diagram.  
\end{theorem}

\med

  \begin{proof} 
 
  We   show that if $\gamma :[0,1] \to\spq$ is a path of chord diagrams, then
$\mu_{\gamma (0)} = \mu_{\gamma (1)}$.   By the connectivity of $\spq$ (theorem \ref{connect}), this will prove the theorem.  To  do this we parameterize
the construction of the operation. Namely, let 
$$Map_*(\gamma , M) = \{(t, f): \, t \in [0,1], \, f \in Map_*(\gamma (t), M) \}.$$
Then there are restriction maps to the incoming and outgoing boundaries,
$\rho_{in} : Map_*(\gamma , M) \to (LM)^p$,  and $\rho_{out} : Map_*(\gamma , M) \to (LM)^q$.
Let $p : Map_*(\gamma , M) \to [0,1]$ be the projection map.  Then lemma \ref{nbd} implies the following.

\begin{lemma}  The product $\rho_{in} \times p : Map_*(\gamma , M) \hk (LM)^p \times [0,1]$ is a codimension $-\chi(\sig_c)n$  embedding, and has an open
neighborhood
 $\nu (\gamma)$  diffeomorphic to the total space of
the  vector bundle whose fiber  over $(t,f) \in Map_*(\gamma , M) $ is
 given by
$$
\nu (\gamma)_{(t,f)} = \bigoplus_{v \in \sigma (\gamma (t))} \bigoplus_{(\mu (v) -1)} T_{f(v)}M.
$$
\end{lemma}
 This allows us to define a Thom collapse map, $\tau : (LM)^p \times [0,1] \to (Map_*(\gamma , M))^{\nu (\gamma)}$
which defines a homotopy between the collapse maps $\tau_0 : (LM)^p \to Map_*(\gamma (0), M)^{\nu  (\gamma (0))} \hk Map_*(\gamma , M)^{\nu 
(\gamma  )}$ and  $\tau_1 : (LM)^p \to Map_*(\gamma (1), M)^{\nu  (\gamma (1))} \hk Map_*(\gamma , M)^{\nu  (\gamma ) }$.

One can then define the push-forward map,
$$
\begin{CD}
\rn : h_*((LM)^p \times [0,1]) @>\tau_* >> h_*( (Map_*(\gamma , M))^{\nu (\gamma)})  @>t>\cong> h_{*+\chi \cdot n}( (Map_*(\gamma , M))  
\end{CD}
$$ 
and then an operation
$$
 \mu_\gamma  = (\rho_{out})_* \circ \rn:   h_*((LM)^p \times [0,1])   \to h_{*+\chi \cdot n}( (Map_*(\gamma , M))  \to h_{*+\chi \cdot n} ((LM)^q).
$$
The restriction  of this  operation to $h_*((LM)^p \times \{0\}) \hk h_*((LM)^p \times [0,1]) $ is, by definition, $\mu_{\gamma (0)}$, and the restriction
to $h_*((LM)^p \times \{1\}) \hk h_*((LM)^p \times [0,1]) $ is $\mu_{\gamma (1)}$.  This proves that these two operations are equal.   \end{proof}

 \med

 Now that we have theorem \ref{invariant} we can  introduce the  notation $\mu_{p,q}(g)$ to stand for $\mu_c : h_*(LM)^{\otimes p} \to h_*(LM)^{\otimes
q}$ for any Sullivan chord diagram $c \in \spq$.  $\mu_{p,q}(g)$ is an operation that lowers the total degree by $(2g-2+p+q)n$.  

\med
\bf Remark. \rm  The above argument is easily modified to show that any element $\alpha \in h_*(\spq) $ defines a string topology operation
$\mu_{p,q}(g)(\alpha)$.   The operations we are dealing with correspond to the class $1\in h_0(\spq)$.   

\med

In order to complete the proof of theorem \ref{one}, by the  the correspondence between two dimensional TQFT's and Frobenius algebras 
\cite{dijkgraaf},
\cite{abrams}, it suffices to show that these operations respect the gluing of surfaces.   

\begin{theorem} \label{compose}
 $\mu_{q,r}(g_2) \circ \mu_{p,q}(g_1) = \mu_{p,r}(g_1+g_2+q-1) : h_*(LM)^{\otimes p} \to h_*(LM)^{\otimes q} \to h_*(LM)^{\otimes r}.$
 
\end{theorem}
\begin{proof}  Let $c_1 \in  \cc \cf^\mu_{p,q}(g_1)$ and $c_2 \in  \cc\cf^\mu_{q,r}(g_2)$.  Notice that we can glue $c_1$ to $c_2$ to obtain a Sullivan chord
diagram in $c_1\#c_2 \in \cc \cf^\mu_{p,r}(g_1 +g_2+q-1)$ in the following way.   

Identify  the outgoing boundary circles of $c_1$ with the incoming
boundary circles of $c_2$ using the parameterizations,  and input the vertices and ghost edges of $c_2$ into the
diagram $c_1$ using these identifications.   Figure \ref{gluing} gives an example of this gluing procedure with 
$c_1 \in  \cc \cf^\mu_{1,2}(0)$, $c_2 \in \cc \cf^\mu_{2,2}(0)$, and $c_1 \# c_2 \in \cc \cf^\mu_{1,2}(1)$. For clarity the vertices have been labeled in
these diagrams, both before and after gluing.

\begin{figure}
\begin{center}
\mbox{
\subfigure[$c_1$]{\epsfig{file=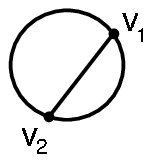, height=60pt, width=60pt}}\qquad
\subfigure[$c_2$]{\epsfig{file=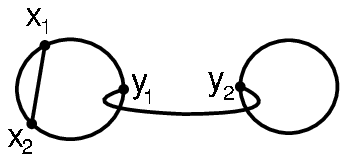, height=60pt, width=120pt}}\qquad
\subfigure[$c_1 \# c_2$]{\epsfig{file=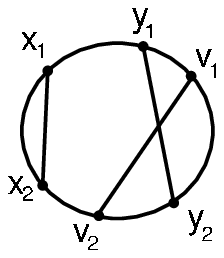, height=60pt, width=60pt}}}
\caption{Gluing $c_1$ and $c_2$} \label{gluing}
\end{center}
\end{figure}

\bf Note.  \sl We are not claiming that this gluing procedure is continuous, or even well-defined.  The ambiguity in definition occurs
if, when one identifies the outgoing boundary circle of $c_1$ with an  incoming boundary circle of $c_2$, a circular vertex $x$ of
$c_2$ coincides with  a circular vertex $v$ of $c_1$ that lies on a ghost edge in the boundary cycle.   Then there is an
ambiguity over whether to place $x$ at $v$ or at the other vertex of the ghost edge.  However for our purposes, we can make any
such choice, since the operations that two such glued surfaces define are equal, by theorem \ref{invariant}.

\rm
Notice that the parameterizations give us maps of the collapsed fat graphs,
$$\phi_1 :S( c_1) \to S( c_1 \#c_2) \quad  {\rm and} \quad \phi_2 : S(c_2) \to S( c_1\#c_2).$$
These induce  a diagram of mapping spaces,
$$
\begin{CD}
Map_*(c_2, M)@<\bph_2 << Map_*(c_1\#c_2, M) @>\bph_1 >> Map_*(c_1, M)
\end{CD}
$$

The next two lemmas  follow from a verification of the definitions of the mapping spaces and the maps $\phi_i$.

\med
\begin{lemma}\label{emb}
$\bph_1 : Map_*(c_1\#c_2, M) \to Map_*(c_1, M)$ is an embedding, whose image has a neighborhood diffeomorphic
to the total space of the bundle $\bph_2^*(\nu (c_2))$, where $\nu(c_2)\to Map_*(c_2, M) $ is the normal bundle of
$\rho_{in} : Map_*(c_2, M)  \hk (LM)^p$ described in lemma \ref{nbd}.
\end{lemma}

This allows the definition of a Thom collapse map $\tau_{\phi_1} : Map_*(c_1, M) \to Map_*(c_1\#c_2, M)^{ \bph_2^*(\nu (c_2))} $
and therefore a push-forward map in homology,
$(\bph_1)_! : h_*(Map_*(c_1, M)) \to
h_{*+ \chi (c_2)n} (Map_*(c_1\#c_2, M))$.

\med
\begin{lemma}\label{commute}  The following diagram commutes:
$$
\begin{CD}
(LM)^r   @>=>> (LM)^r \\
@A\rho_{out}(2) AA   @AA\rho_{out}(1\#2)A \\
Map_*(c_2, M)  @<\bph_2 <<  Map_*(c_1 \#c_2, M)  @>\rho_{in}(1\#2) >>  (LM)^p\\
 @V\rho_{in}(2) VV  @VV\bph_1 V  @VV = V \\
 (LM)^q  @<<\rho_{out}(1) <   Map_*(c_1, M)   @>>\rho_{in}(1)>  (LM)^p.
 \end{CD}
$$ The indexing of the restriction maps corresponds to the indexing of the chord diagrams in the obvious way.
\end{lemma}

By the naturality of the Thom collapse map, and therefore the homological pushout  construction, we therefore
have the following corollary.

\begin{corollary}

$ 1.  \,  (\rho_{in} (1\#2))_! = (\bph_1)_! \circ (\rho_{in} (1))_! : h_*((LM)^p) \to h_{*+\chi (c_1\#c_2)\cdot n}(Map_*(c_1\#c_2,
M))$

$2. \,  (\bph_2)_* \circ (\bph_1)_!  = (\rho_{in}(2))_! \circ (\rho_{out}(1))_* : h_*(Map(c_1, M)) \to h_{*+\chi (c_2)\cdot
n}(Map_*(c_2, M)) $

$3.   (\rho_{out}(1\#2))_*  = (\rho_{out}(2))_* \circ (\bph_2)_* : h_*(Map_*(c_1\#c_2, M)) \to h_*((LM)^r).$
 \end{corollary}

We may now complete the proof of theorem \ref{compose}.  We have

\begin{align}
\mu_{p,r}(g_1+g_2+q-1) = \mu_{c_1\#c_2} &= (\rho_{out}(1\#2))_* \circ (\rho_{in}(1\#2))_!  \notag\\
&= (\rho_{out}(1\#2))_*  \circ  (\bph_1)_! \circ (\rho_{in} (1))_!  \notag\\
&=  (\rho_{out}(2))_* \circ (\bph_2)_*\circ (\bph_1)_! \circ (\rho_{in} (1))_!  \notag \\
&= (\rho_{out}(2))_* \circ  (\rho_{in}(2))_! \circ (\rho_{out}(1))_* \circ (\rho_{in} (1))_!  \notag \\
&= \mu_{c_2}\circ \mu_{c_1} \notag  \\
&= \mu_{q,r}(g_2) \circ \mu_{p,q}(g_1).  \notag
\end{align}

\end{proof}
As  observed in \cite{abrams}, a Frobenius algebra without counit is the same thing as a  positive boundary 
topological quantum field
theory.  We  have now verified that the string topology operations define such a theory for any generalized cohomology
theory $h^*$ satisfying the conditions described above.  Recall that it was observed in   \cite{chassull} \cite{cj}, that the unit in the algebra structure of
$h_*(M)$ is the fundamental class, $[M] \in h_n(M) \hk h_n(LM)$, where the second map is induced by the inclusion of the manifold in the loop space as the
constant loops,
$\iota : M \hk LM$.   Thus $h_*(LM)$ is a unital Frobenius algebra without a counit.   This proves theorem \ref{one}.

\med

\section{Capping off boundary components: issues surrounding the unit and counit }

 The unit in the Frobenius algebra stucture can be constructed in the same way as 
 the other string topology operations as follows.

Consider the disk $D^2$ as a surface with zero incoming boundary component and one outgoing boundary component.  A graph
$c_D$ that represents $D^2$ can be taken to be a point (i.e a single vertex).  Formally, the restriction to the zero incoming
boundary components is the map
$$\rho_{in} : M = Map_*(c_D, M) \to Map(\emptyset , M) = point.
$$
The push-forward map in this setting
$$
\rn : h_*(point) \to h_{*+n}(M)
$$
is the $h_*$-module map defined by sending the generator to $[M] \in h_n(M)$.
The restriction to the outgoing boundary component  is the map
$$\rho_{out} : M = Map_*(c_D, M) \to LM$$ which is given by $\iota : M \hk LM$.  Thus the unit is given by
the
$h_*$ module homomorphism
$$
\mu_{D^2} = (\rho_{out})_* \circ \rn = \iota_* \circ \rn : h_* \to h_n(M) \to h_n(LM)
$$
which sends the generator to the fundamental class.  

\med
The issue of the existence (or nonexistence) of a \sl counit \rm in the Frobenius algebra structure given by theorem \ref{one}
is formally the same (or dual) to the existence of a unit, but is geometrically much more difficult and subtle.   Namely, for this
operation one must consider $D^2$ as a surface with one incoming boundary, and zero outgoing boundary components.  In this
setting the roles of the restriction maps $\rho_{in}$ and $\rho_{out}$ are reversed , and one obtains the diagram
$$
\begin{CD}
 Map(\emptyset, M)  @<\rho_{out} <<  Map_*(c_D, M) @>\rho_{in} >> LM
\end{CD}
$$
or, equivalently,
$$
\begin{CD}
point  @<\epsilon <<  M @>\iota >> LM.
\end{CD}
$$where $\epsilon : M \to point$ is the constant map.

Now notice that in this case, unlike when any of the other fat graphs were considered, the embedding of
$Map_*(c_D, M) \hk LM$  (i.e $\iota : M \hk  LM$) is of \sl infinite \rm codimension.  Therefore to find a theory $h_*$ that supports
a counit in the Frobenius algebra structure of $h_*(LM)$, one needs to be able to define a push-forward map for this infinite
codimensional embedding.   Now   in their work on genera of loop spaces, \cite{andomorava},  Ando and Morava argued that if one has   
a theory where this push-forward map exists,  one would  need that  the Euler class of the normal bundle
$e(\nu (\iota))
\in h^*(M)$ is invertible.  So let us now consider this normal bundle. 

The embedding of $M$ as the constant loops in $LM$ is $S^1$ -equivariant where $S^1$ acts triviallly
on $M$. When $M$ is a simply connected almost complex manifold, the normal bundle has the following
description (see \cite{andomorava}, for example). 
 \med
\begin{lemma}\label{normal}  The normal bundle $\nu (\iota) \to M$ of the embedding $\iota : M \hk LM$ is equivariantly isomorphic
to the direct sum,
$$
\nu (\iota) \cong \bigoplus_{k \neq 0} TM \otimes_{\bc} \bc (k)
$$
where $\bc (k)$ is the one dimensional representation of $S^1$ of weight $k$.
\end{lemma}

This says that  the Euler class of the normal bundle  will have the formal description
\begin{equation}\label{euler}
e(\nu(\iota)) = \prod_{k\neq 0} e(TM \otimes \bc (k)).
\end{equation}
Thus a theory $h_*(LM)$ that supports a counit in a Frobenius algebra structure   should  have the following properties.
\begin{enumerate}
\item $h_*$ should be an $S^1$ - equivariant theory in order to take advantage of the different equivariant structures
of the summands $TM \otimes \bc (k)$.
\item $h_*$ should be a ``pro-object" - an inverse system of homology groups, so that it can accommodate this infinite product.
\end{enumerate}

\section{The polarized loop space and its Atiyah dual}

Motivated by these homological requirements, in this section we show that
 the loop space of an almost complex manifold  has a natural equivariant pro-object (a
``prospectrum") associated to it.   The ideas for the constructions in this section stem from
conversations with Graeme Segal.
 Throughout this section we assume that $M$ is a simply connected, oriented, closed $n$-manifold.

Let $-TM$ be the virtual bundle ($K$-theory class) given by the opposite of the tangent bundle.  Let $\mtm$ be  its Thom
 spectrum.  We refer to   $\mtm$ as the  ``Atiyah dual" of $M_+$ because of Atiyah's well known theorem stating
that $\mtm$ is equivalent to the Spanier Whitehead dual of $M_+$.  (Here $M_+$ is $M$ together with a disjoint
basepoint.)  This gives
$\mtm$ the structure of a ring spectrum,  whose multiplication $m : \mtm \wedge \mtm \to \mtm$ is dual to the diagonal
$\Delta : M
\to M\times M$.  When one applies homology and the Thom isomorphism, this multiplication realizes the intersection product
($\cap$), meaning that the following diagram commutes.
$$
\begin{CD}
H_*(\mtm) \otimes H_*(\mtm)  @>m_* >> H_*(\mtm) \\
@Vt V\cong V   @V\cong Vt V \\
H_{*+n} (M) \otimes H_{*+n}(M)  @>>\cap >  H_{*+n} (M).
\end{CD}
$$
Here $H_*(\mtm)$ refers to the spectrum homology of $\mtm$.

It is therefore natural to expect  that an appropriate pro-object that carries the string topology operations,  including a counit,
(i.e a  2 dimensional TQFT, or Frobenius algebra structure), would be a prospectrum model for the Atiyah dual of the loop space, 
$LM^{-TLM}$.  

In studying homotopy theoretic aspects of symplectic Floer homology, the first author, Jones and Segal used   pro-spectra
associated to certain infinite dimensional bundles \cite{cjs}.   The construction was the following.  If $E \to X$ is an infinite dimensional
bundle with a filtration by finite dimensional subbundles,
$$
\cdots \hk E_i \hk E_{i+1}\hk \cdots E
$$
such that $\bigcup_i E_i$ is a dense subbundle of $E$, then one can define the prospectrum $X^{-E}$ to be the inverse system,
$$
\begin{CD}
\cdots \leftarrow X^{-E_{i-1}} @<u_{i}<<  X^{-E_i} @<u_{i+1}<<  X^{-E_{i+1}} \leftarrow \cdots
\end{CD}
$$
where  $u_j : X^{-E_j} \to X^{-E_{j-1}}$ is the map defined as follows.  Let $e_j :E_{j-1} \hk E_j$ be the inclusion.  Assume for
simplicity that $E_j$ is embedded in a large dimensional trivial bundle, and let $E_j^\perp$ and $E_{j-1}^\perp$ be the
corresponding
orthogonal complements.  One then has an induced  inclusion of complements, $e_j^\perp : E_{j }^\perp \to E_{j-1}^\perp$.  The
induced map of Thom spaces then defines a map of Thom spectra,   
$u_j: X^{-E_j} \to X^{-E_{j-1}}$. A standard homotopy theoretic technique allows one to define this map of Thom spectra even
if $E_j$ is not embeddable in a trivial bundle, by restricting  $E_j$ to finite subcomplexes of $X$ where it is.   

Under the assumption that $M$ is an almost complex  manifold of dimension $n = 2m$, then we are dealing with
 an infinite dimensional vector bundle ($TLM$) whose structure group is the loop group $LU(m)$.  In \cite{cs}, Cohen and Stacey
studied obstructions to finding an appropriate filtration  of an infinite dimensional $LU(m)$ bundle.  In particular for $TLM \to LM$
it was proved that if such a filtration (called a ``Fourier decomposition" in \cite{cs}) exists, then the holonomy of any  unitary
connection on $TM$, 
$h: \Omega M \to U(m)$ is null homotopic.  This ``homotopy flat" condition is far too restrictive for our purposes, but we can
get around this problem by taking into account the canonical \sl polarization \rm of the tangent bundle  $TLM$ of an almost complex
manifold. 

\med
Recall  that a  \sl polarization \rm of a Hilbert space $E$ is an equivalence  class of decomposition, $E = E_+ \oplus E_-$,
where two such decompositions $E_+ \oplus E_-  = E^\prime_+     \oplus E^\prime_-$ are equivalent if the composition
$E_+ \hk E \to E^\prime_+$ is Fredholm, and $E_+ \hk E \to E^\prime_-$ is compact.  (See \cite{presseg} for details.)
  The restricted general linear group of a polarized space $GL_{res}(E)$ consists of all elements of $GL(E)$
that preserve the polarization.

A polarized vector bundle $\zeta \to X$ is one where every fiber is polarized, and the structure group reduces to the restricted
general linear group.  If $M^{2m}$ is an almost complex manifold, and $\gamma \in LM$, then the tangent space  
$T_\gamma LM$   is the space of $L^2$ vector fields of  $M$ along $\gamma$, and the operator
$$
j\frac{d}{d\theta} : T_\gamma LM \to T_\gamma LM
$$
is a self adjoint Fredholm operator.   Here $\frac{d}{d\theta}$ is the covariant derivative, and $j$ is the almost complex structure.
The spectral decomposition of  $j\frac{d}{d\theta}$ polarizes the bundle $TLM$ according to its positive and negative eigenspaces.  The structure group in
this case is $GL_{res}(L^2(S^1, \bc^m))$, where the loop space $L^2(S^1, \bc^m)$ is polarized according to the Fourier
decompostion. That is, we write
$$
L^2(S^1, \bc^n) = H_+ \oplus H_-  
$$
where $H_+ = Hol (D^2, \bc^n)$  is the space of holomorphic maps of the disk, and $H_-$ is the orthogonal complement.

For a polarized space $E$,  recall from \cite{presseg} that the \sl  restricted  Grassmannian \rm $Gr_{res}(E)$ consists of closed
subspaces $W \subset E$ such that the projections $W\hk E \to E_+$ is Fredholm, and $W\hk E \to E_-$ is Hilbert-Schmidt.
  In the case
under consideration, the  tangent space $ T_\gamma LM$,  is a $L\bc$- module, and therefore a module over the Laurent
polynomial ring, $\bc [z, z^{-1}]$.  Define $\gro (T_\gamma LM) \subset Gr_{res}(T_\gamma LM)$ to be the subspace
 
$$
\gro(T_\gamma LM) = \{W \in Gr_{res}(T_\gamma LM) \, :  \, zW \subset W\}.
$$

\med
For $M^{2m}$ a simply connected almost complex manifold, we can  then define the    \emph{polarized loop space}  
$\lpm$ to be the space 
\begin{equation}\label{polar}
\lpm = \{(\gamma, W) \, : \,  \gamma \in LM,  \ W \in Gr^0_{res}(T_\gamma LM) \}.
\end{equation}

\med
We now consider the $S^1$-equivariance properties of $\lpm$.   The following theorem will be an easy consequence of the results of   
\cite{presseg}, chapter  8.  
 
\begin{theorem} The natural projection  $p : \lpm \to L(M^{2m})$  is an $S^1$ equivariant fiber bundle with fiber
diffeomorphic to the based loop space, $\Omega U(m)$.  The $S^1$-fixed points of $\lpm$  form a bundle over $M$
with fiber the space of group homomorphisms, $Hom (S^1, U(m))$.
 \end{theorem}

\begin{proof} Let $\gamma \in LM$.  The tangent space, $T_\gamma LM = \Gamma_{S^1}(\gamma^*(TM))$, is the space of $L^2$
sections of the pullback of the tangent bundle over the circle.  The $S^1$- action on $LM$ differentiates to make the
tangent bundle $TLM$ an  $S^1$-equivariant bundle.  If $\sigma
\in T_\gamma LM$, and $t\in S^1$ then $t \sigma \in T_{t \gamma} LM$ is defined by $t \sigma (s) = \sigma (t+s)$.
 Since this action preserves the polarization, it induces an action  
\begin{align}
S^1 \times \lpm  &\to \lpm \\
t \times (\gamma, W) &\to (t  \gamma, t W) 
\end{align}
where $tW = \{t\sigma \in T_{t\gamma}LM \, : \, \sigma \in W \}$.   The fact that the projection map
$p : \lpm \to LM$ is an $S^1$-equivariant bundle is clear.   The fiber of this bundle can be identified with
$Gr_{res}^0(L^2(S^1, \bc^m))$, which was proved in \cite{presseg} to be diffeomorphic to $\Omega U(m)$. 
The induced action on $\Omega U(m)$ was seen in \cite{presseg}, (chapter 8) to be given as follows.  For $t \in 
S^1 = \br/\bz$, 
and $\omega \in \Omega U(m)$, $t\cdot \omega (s) = \omega (s+t)\omega(t)^{-1}.$  The fixed points of this
action are the group homomorphisms, $Hom (S^1, U(m))$.  The theorem now follows. \end{proof}

\med
\bf Remark.  \rm  Since the group homomorphisms, $Hom (S^1, U(m))$ are well understood,  one can view the above theorem as saying that the equivariant
homotopy type of $\lpm$ is directly computable in terms of the equivariant homotopy type of $LM$.   

\med
By this theorem, the pullback of the tangent bundle, $p^*TLM \to \lpm$ is an $S^1$ equivariant bundle.  Our final result 
implies   that even though one cannot generally  find a prospectrum modeling the Atiyah dual $LM^{-TLM}$, one can find a
pro-spectrum model of the ``polarized Atiyah dual", $\lpm^{-TLM}$. 

The following theorem says that one can build up  the bundle $p^*(TLM)\to \lpm$ by finite dimensional subbundles.

\med
\begin{theorem}  There is a doubly graded collection of finite dimensional, $S^1$-equivariant
subbundles  of $p^*(TLM) \to \lpm$,
$$
E_{i, j} \to \lpm, \quad i < j
$$
satisfying the following properties:
\begin{enumerate}
\item There are inclusions of subbundles 
$$
E_{i,j} \hk E_{i-1, j}  \quad {\rm and} \quad E_{i,j} \hk E_{i, j+1}
$$
such that $\bigcup_{i,j} E_{i,j}$ is a dense subbundle of $p^*(TLM)$. 
\item  The subquotients,
$$
 E_{i-1, j} /E_{i,j}  \quad {\rm and} \quad E_{i, j+1}/ E_{i,j}
$$
are $m$ dimensional $S^1$-equivariant complex vector bundles that are
nonequivariantly isomorphic to the pullback of the  tangent bundle $\tilde p^* TM$,
where $\tilde p : \lpm \to M$ is the composition of $p : \lpm \to LM$ with the map
$e_1: LM \to M$ that evaluates a loop at the basepoint $1 \in S^1$. 
\end{enumerate}
\end{theorem}

\med
{\bf Remark.}  Such a filtration is a ``Fourier decomposition" of the loop bundle $p^*(TLM) \to \lpm$
as defined in \cite{cs}
\begin{proof}
We first define certain infinite dimensional subbundles $E_i \subset p^*(TLM) \to \lpm$.    Define 
the fiber over $(\gamma, W) \in \lpm$ to be
$$
(E_i)_{(\gamma, W)} = z^{-i}W \subset T_{\gamma}(LM).
$$
We note that $E_i$ is an equivariant subbundle,  with the property that $zE_i \subset E_i$.  Furthermore there is a filtration of subbundles
$$
\cdots \hk E_i \hk E_{i+1} \hk \cdots p^*(TLM)
$$
with $\bigcup_i E_i$ a dense subbundle of $p^*(TLM)$.  Notice that for $j > i$,  the subquotient $E_{j}/E_{ i}$ has fiber
at $(\gamma, W)$ given by
$z^{-j }W \cap (z^{-i}W)^\perp$  where $ (z^{ -i}W)^\perp \subset T_\gamma LM$ is the orthogonal complement of
$z^{ -i }W$.   For $j-i =
 1$, an easy argument (done in \cite{cs})  gives that the composition 
\begin{equation}\label{eval}
\begin{CD}
z^{- j  }W \cap (z^{-(j-1)}W)^\perp \hk T_\gamma LM  @>e_1 >>T_{\gamma (1)}(M)
\end{CD}
\end{equation}
is an isomorphism.   
For $i < j$ we define the bundle
$E_{i,j} \to \lpm$ to be the quotient $E_j/E_{ i}$.  It is the vector bundle whose fiber over $(\gamma, W)$ is $z^{-j }W
\cap (z^{- i }W)^\perp$. By (\ref{eval}), the subquotient of the bundle $E_{j-1.j}$ is (nonequivariantly)  isomorphic to the
pullback of the tangent bundle $TM
\to M$  under the composition $\begin{CD}\lpm  @>p >> LM  @>e_1 >> M\end{CD}$.   In general the bundle $E_{i,j}$ is nonequivariantly isomorphic to
the Whitney sum of $j-i$ copies of $ \tilde p^*(TM)$. 

Now since $z^{-j}W \cap (z^{-i}W)^\perp $ is a subspace  of both $z^{-(j+1)}W \cap (z^{-i}W)^\perp $ and of $z^{-j}W \cap
(z^{-(i-1)}W)^\perp$, we have inclusions    $E_{i,j} \hk E_{i, j+1}  \quad {\rm and} \quad E_{i,j} \hk E_{i-1, j}.$  Clearly
$\bigcup_{i,j} E_{i,j}$ is a dense subbundle of
$p^*(TLM)$.  The theorem follows.
   \end{proof}
  
\med

Since the bundles $E_{i,j} \to \lpm$ are finite dimensional $S^1$-equivariant bundles,  we can construct the Thom spectrum of the $S^1$-equivariant
virtual bundle,  $-E_{i,j}$, which we denote by $(\lpm)^{-E_{i,j}}$.  Notice  the inclusions of bundles $E_{i,j} \hk E_{i,j+1}$ and 
$E_{i,j} \hk E_{i-1,j}$ induce maps of virtual bundles, $\tau_{i,j} : - E_{i,j+1} \to -E_{i,j}$ and $\sigma_{i,j} : -E_{i-1,j} \to
-E_{i,j}$, which yields an inverse system of  
$S^1$-equivariant spectra,
$$
\begin{CD}
\vdots  && \vdots \\
@V \sigma_{i-1,j} VV  @VV\sigma_{i-1,j+1} V  \\
(\lpm)^{-E_{i-1,j}}   @<\tau_{i-1,j}<< (\lpm)^{-E_{i-1,j+1}} @<\tau_{i-1, j+1}<< \cdots  \\
@V \sigma_{i,j} VV  @VV\sigma_{i,j+1} V  \\
(\lpm)^{-E_{i,j}}   @<\tau_{i,j}<< (\lpm)^{-E_{i,j+1}} @<\tau{i, j+1}<< \cdots \\
@V \sigma_{i+1,j} VV  @VV\sigma_{i+1,j+1} V  \\
\vdots  && \vdots 
\end{CD}
$$
This system defines a pro-object in the category of $S^1$-equivariant spectra that we call the polarized Atiyah dual,
$\lpm^{-TLM}$.  If one applies an equivariant  homology theory   to this prospectrum, one
gets a pro-object in the category of graded abelian groups.  Notice that in cohomology, the structure maps $\tau_{i,j}$ and
$\sigma_{i,j}$ will induce multiplication by the equivariant Euler classes of the orthogonal complement bundles
of these inclusions.  As seen above, these orthogonal complement bundles are nonequivariantly isomorphic to the pull back
of the tangent bundle, $TM$.  However they have different equivariant structures.   In future work we will study
those equivariant cohomology theories for which these Euler classes are units, with the goal being to prove that
  such theories, when applied to this prospectrum,   support the string topology operations including one that corresponds to a
disk viewed as a cobordism from a circle to the empty set.  By gluing, this will allow the construction of string topology
operations for  closed surfaces as well as surfaces with a positive number of outgoing boundary components.

\end{document}